\documentclass[11pt]{article} 
\usepackage{amssymb,amsmath} 

\def\R{\mathbb{R}}
\def\Z{\mathbb{Z}}
\def\B{{\cal B}}
\def\L{{\cal L}}
\def\P{{\cal P}}
\def\T{{\cal T}}
\def\A{{\bf A}} 
\def\a{{\bf a}}
\def\b{{\bf b}}
\def\c{{\bf c}}

\def\t{{\bf t}} 
\def\v{{\bf v}} 
\def\w{{\bf w}} 
\def\x{{\bf x}} 
\def\z{{\bf z}} 
\def\1{{\bf 1}} 
\newcommand\vol{\operatorname{vol}}

\begin{document}
\setlength{\parindent}{0pt}
\setlength{\parskip}{0.4cm}
\bibliographystyle{amsplain} 

\newtheorem{theorem}{Theorem}
\newtheorem{proposition}[theorem]{Proposition}
\newtheorem{corollary}[theorem]{Corollary}

\begin{center}

\Large{\bf The Ehrhart polynomial of the Birkhoff polytope} 
\footnote{Appeared in {\it Discrete \& Computational Geometry} {\bf 30}, no.~4 (2003), 623--637.} 
\normalsize

{\sc Matthias Beck and Dennis Pixton} 

\end{center}

\small
\begin{quote}
{\it All means (even continuous) sanctify the discrete end.} \ 
Doron Zeilberger\footnote{E-mail to the first author on {\tt Wed, 14 Aug 2002 16:29:35 -0400 (EDT)}} 
\end{quote}
\normalsize

\footnotesize
\emph{Abstract}: The $n^{\text{th}}$ Birkhoff polytope is the set of all doubly stochastic 
$n \times n$ matrices, that is, those matrices with nonnegative real coefficients in which 
every row and column sums to one. A wide open problem concerns the volumes of these polytopes, 
which have been known for $n \leq 8$. We present a new, complex-analytic way to compute the 
Ehrhart polynomial of the Birkhoff polytope, that is, the function counting the integer points in the 
dilated polytope. One reason to be interested in this counting function is that the 
leading term of the Ehrhart polynomial is---up to a trivial factor---the volume of the polytope. 
We implemented our methods in the form of a computer program, which yielded the Ehrhart polynomial 
(and hence the volume) of the ninth and the volume of the tenth Birkhoff polytope.
\normalsize


\section{Introduction} 

One of the most intriguing objects of combinatorial geometry is the \emph{$n^{\text{th}}$ 
Birkhoff polytope} 
  \[ \B_n = \left\{ \left( \begin{array}{ccc} x_{11} & \cdots & x_{1n} \\ \vdots & & \vdots \\ x_{n1} & \dots & x_{nn} \end{array} \right) \in \R^{n^2} : \ x_{jk} \geq 0 , \begin{array}{l} \sum_j x_{jk} = 1 \text{ for all } 1 \leq k \leq n \\ \sum_k x_{jk} = 1 \text{ for all } 1 \leq j \leq n \end{array} \right\} \ , \] 
often described as the set of all $n \times n$ \emph{doubly stochastic matrices}. 
$\B_n$ is a convex polytope with integer vertices. It possesses fascinating combinatorial 
properties \cite{billerasarangarajan,brualdigibson1,brualdigibson3,chanrobbinscatalan,zeilbergercatalan} and relates to many mathematical areas \cite{diaconisgangoli,knuthmagic}. 
A long-standing open problem is the determination of the relative volume 
of $\B_n$, which had been known only up to $n=8$ \cite{chanrobbinsbirkhoff,mount}. 
In this paper, we propose a new method of calculating this volume and use it to compute 
$\vol \B_9$ and $\vol \B_{10}$. 

One of the recent attempts to compute $\vol \B_n$ relies on the theory of counting 
functions for the integer points in polytopes. Ehrhart proved \cite{ehrhart2} that for 
a polytope $\P \subset \R^d$ with integral vertices, the number 
  \[ L_\P(t) := \# \left( t \P \cap \Z^d \right) \] 
is a polynomial in the positive integer variable $t$. 
He showed various other properties for this counting function (in fact, in the more general setting 
of $\P$ having rational vertices), of which we mention three here: 
\begin{itemize} 
  \item The degree of $L_\P$ is the dimension of $\P$. 
  \item The leading term of $L_\P$ is the relative volume of $\P$, normalized with respect to the sublattice of $\Z^d$ on the affine subspace spanned by $\P$. 
  \item Since $L_\P$ is a polynomial, we can evaluate it at nonpositive integers. These evaluations yield 
    \begin{align} &L_\P(0) = \chi (\P) \ , \label{constant} \\ 
                  &L_\P(-t) = (-1)^{ \dim \P } L_{\P^\circ}(t) \ . \label{rec} \end{align} 
\end{itemize} 
Here $\chi (\P)$ denotes the Euler characteristic, $\P^\circ$ the relative interior of $\P$. 
The reciprocity law (\ref{rec}) was in its full generality proved by Macdonald \cite{macdonald}. 

The application of this theory to the Birkhoff polytope $\B_n$ incorporates the nice interpretation 
of the number of integral points in $ t \B_n $ as the number of \emph{semi-magic squares}, namely, 
square matrices whose nonnegative integral coefficients sum up to the same integer $t$ along 
each row and column. 

We will denote the Ehrhart polynomial of $\B_n$ by 
  \[ H_n (t) := L_{\B_n} (t) \ . \] 
It is not hard to see that $\dim \B_n = (n-1)^2$, hence $H_n$ is a polynomial in $t$ of degree 
$(n-1)^2$. The first two of these polynomials are trivial: 
  \[ H_1 (t) = 1 \ , \qquad \qquad H_2(t) = t+1 \ , \] 
the first nontrivial case was computed by MacMahon \cite{macmahon} as 
  \begin{equation}\label{h3} H_3 (t) = 3 {{ t+3 }\choose{ 4 }} + {{ t+2 }\choose{ 2 }} \ . \end{equation} 
The structural properties of $H_n$ were first studied in \cite{ehrhartmagic,stanleymagic,stein}. 
It is a nice exercise to deduce from (\ref{rec}) that 
  \begin{equation}\label{reczeros} H_n (-n-t) = (-1)^{ n - 1 } \ H_n (t) \end{equation} 
and 
  \[ H_n (-1) = H_n (-2) = \dots = H_n (-n+1) = 0 \ . \] 
This allows the following strategy of computing $H_n$, and therefore, the volume of $\B_n$: 
compute the first $\binom {n-1} 2$ values 
of $H_n$, use the above symmetry and trivial values of $H_n$, and calculate the polynomial $H_n$ 
by interpolation. In fact, as far as we are aware of, the volume of $\B_8$ was 
computed using essentially this method, combined with some nice computational tricks \cite{chanrobbinsbirkhoff,mount}. 

We propose a new, completely different approach of computing $H_n$ (and hence $\vol \B_n$). 
It is based on an analytic method by the first author of computing the Ehrhart polynomial of a 
polytope \cite{tetrahed}. We will introduce the application of this method to the Birkhoff polytope 
in the following section. 

Some recent refreshing approaches---of a more algebraic-geometric/topological flavor---to the problem of computing $\vol \B_n$ and $H_n$ 
can be found in \cite{baldonivergne,deloerasturmfels,sturmfelsmagic}. 


\section{An integral counting integers} 

A convex polytope $\P \subset \R^d$ is an intersection of halfspaces. This allows the 
compact description 
  \[ \P = \left\{ \x \in \R^d : \ \A \x \leq \b \right\} \ , \] 
for some $(m \times d)$-matrix $\A$ and $m$-dimensional vector $\b$. Here the inequality is 
understood componentwise. In fact, we may convert these inequalities into equalities by 
introducing `slack variables.' If $\P$ has rational vertices (those polytopes are called 
\emph{rational}), we can choose $\A$ and 
$\b$ in such a way that all their entries are integers. In summary, we may assume that a 
convex rational polytope $\P$ is given by 
  \begin{equation}\label{polytope} \P = \left\{ \x \in \R_{ \geq 0 }^d : \ \A \x = \b \right\} \ , \end{equation} 
where $\A \in M_{m \times d} (\Z)$ and $\b \in \Z^m$. 
(If we are interested in counting the integer points in $\P$, we may assume that 
$\P$ is in the nonnegative orthant, i.e., the points in $\P$ have nonnegative coordinates, 
as translation by an integer vector does not change the lattice-point count.) 
The following straightforward theorem can be found in \cite{tetrahed}. 
We use the standard multivariate notation $ \v^{\w} := v_1^{ w_1 } \cdots v_n^{ w_n } $. 
\begin{theorem}\label{intthm} {\rm \cite[Theorem 8]{tetrahed}} Suppose the convex rational polytope $\P$ 
is given by {\rm (\ref{polytope})}, and denote the columns of $\A$ by $\c_1,\dots,\c_d$. Then 
  \[ L_\P (t) = \frac{ 1 }{ (2 \pi i)^m } \int_{ \left| z_1 \right| = \epsilon_1 } \cdots \int_{ \left| z_m \right| = \epsilon_m } \frac{ z_1^{ -t b_1 - 1 } \cdots z_m^{ -t b_m - 1 } }{ \left( 1 - \z^{ \c_1 } \right) \cdots \left( 1 - \z^{ \c_d } \right) } \ d \z \ . \] 
Here  $ 0 < \epsilon_1, \dots , \epsilon_m < 1 $ are distinct real numbers. 
\end{theorem} 
It should be mentioned that $L_\P$ is in general \emph{not} a polynomial if the vertices of $\P$ are not 
integral, but a quasipolynomial, that is, an expression of the form
  \[ c_{d}(t) t^{d} + \dots + c_{1}(t) t + c_{0}(t) \ , \] 
where $ c_{0}, \dots , c_{d} $ are periodic functions in $t$. 
(See, for example, \cite[Section 4.4]{stanleyec1} for more information about quasipolynomials.) 
Theorem \ref{intthm} applies to these slightly 
more general counting functions; however, in this article we will only deal with polytopes with 
integer vertices, for which $L_\P$ is a polynomial. 

We also note here that Theorem \ref{intthm} can be used to quickly compute by hand formulas for certain classes 
of polytopes (see, for example, \cite{tetrahed}). In this project, we take a slightly different 
approach and use this theorem to efficiently derive formulas with the help of a computer. 

We can view the Birkhoff polytope $\B_n$ as given in the form of (\ref{polytope}), where 
  \[ \A = \left( \begin{array}{cccccccccccccc} 
              1 & \cdots & 1  \\ 
                &    &    & 1 & \cdots & 1 \\[-5.5pt] 
                &    &    &   &   &    &   \ddots & \\ 
                &    &    &   &   & & &    1 & \cdots & 1 \\[-4pt] 
              1 &    &    & 1 &    & &   &  1 \\[-5.5pt] 
	        & \ddots & &  & \ddots & &  \cdots &   & \ddots \\ 
		&    &  1 &   &    & 1  &    &   &    & 1 
          \end{array} \right) \] 
is a $(2n \times n^2)$-matrix and $\b = ( 1, \dots, 1 ) \in \Z^{2n}$. 
Hence Theorem \ref{intthm} gives for this special case 
  \[ H_n (t) = \frac 1 { ( 2 \pi i )^{2n} } \int \cdots \int \frac{ ( z_1 \cdots z_{2n} )^{ -t-1 } }{ ( 1 - z_1 z_{n+1} ) ( 1 - z_1 z_{n+2} ) \cdots ( 1 - z_n z_{2n} ) } \ d \z \ . \] 
Here it is understood that each integral is over a circle with radius $<1$ centered at 0; all appearing radii 
should be different. We can separate, say, the last $n$ variables and obtain 
  \[ H_n (t) = \frac 1 { ( 2 \pi i )^{n} } \int \cdots \int ( z_1 \cdots z_{n} )^{ -t-1 } \left( \frac 1 { 2 \pi i } \int \frac{ z^{ -t-1 } }{ ( 1 - z_1 z ) \cdots ( 1 - z_n z ) } \ dz \right)^n \ dz_n \cdots dz_1 \ . \] 
We may choose the radius of the integration circle of the innermost integral to be smaller 
than the radii of the other integration paths. Then this 
innermost integral is easy to compute: It is equal to the residue at 0 of 
  \[ \frac{ 1 }{ z^{ t+1 } ( 1 - z_1 z ) \cdots ( 1 - z_n z ) } \] 
and, by the residue theorem, equal to the negative of the sum of the residues at $z_1^{-1}, \dots, z_n^{-1}$. 
(Note that here we use the fact that $t>0$.) 
The residues at these simple poles are easily computed: the one, say, at $1/z_1$ can be calculated as 
\begin{align*}
&\lim_{z \to 1/z_1} \left( z - \frac 1 {z_1} \right) \frac{ 1 }{ z^{ t+1 } ( 1 - z_1 z ) \cdots ( 1 - z_n z ) } 
  = \lim_{z \to 1/z_1} \frac{ z - \frac 1 {z_1} }{ 1 - z_1 z } \ \frac{ 1 }{ z^{ t+1 } ( 1 - z_2 z ) \cdots ( 1 - z_n z ) } \\ 
&\qquad = - \frac 1 {z_1} \ \frac{ z_1^{ t+1 } }{ ( 1 - z_2/z_1 ) \cdots ( 1 - z_n/z_1 ) } 
  = - \frac{ z_1^{ t+n-1 } }{ ( z_1 - z_2 ) \cdots ( z_1 - z_n ) } \ . \\ 
\end{align*} 
This yields the starting point for our computations. 
\begin{theorem}\label{Hthm} For any distinct $ 0 < \epsilon_1, \dots, \epsilon_n < 1 $, 
  \[ H_n (t) = \frac 1 { ( 2 \pi i )^{n} } \int_{ \left| z_1 \right| = \epsilon_1 } \cdots \int_{ \left| z_n \right| = \epsilon_n } ( z_1 \cdots z_{n} )^{ -t-1 } \left( \sum_{k=1}^n \frac{ z_k^{ t+n-1 } }{ \prod_{j \not= k} (z_k - z_j) } \right)^n dz_n \cdots dz_1 \ . \] 
\end{theorem} 

\emph{Remark}. It can be proved from the form of the integrand that $H_n$ is indeed 
a polynomial in $t$: To compute the integral, one has to execute a (huge) number of limit 
computations, which yield ``at worst" powers of $t$ (as a consequence of L'Hospital's Rule). 
In fact, one can make this property more apparent by noticing that the expression in parenthesis 
is actually a polynomial; namely 
  \[ H_n (t) = \frac 1 { ( 2 \pi i )^{n} } \int_{ \left| z_1 \right| = \epsilon_1 } \cdots \int_{ \left| z_n \right| = \epsilon_n } ( z_1 \cdots z_{n} )^{ -t-1 } \left( \sum_{ m_1 + \dots + m_n = t } z_1^{m_1} \cdots z_n^{m_n} \right)^n dz_n \cdots dz_1 \ , \] 
where the sum is over all ordered partitions of $t$. This formula can also be proved ``more 
directly" combinatorially.\footnote{The authors thank Sinai Robins and Frank Sottile for 
their help in the proof of this equivalence and its combinatorial interpretation.} 


\section{Small $n$ do not require a computer} 

We will now illustrate the computation of $H_n$ (and hence $\vol \B_n$) by means of Theorem 
\ref{Hthm} for $n=3$ and $4$. These calculations ``by hand" give an idea what computational 
tricks one might use in tackling larger $n$ with the aid of a computer. 

By the theorem, 
  \[ H_3(t) = \frac{ 1 }{ (2 \pi i)^3 } \int (z_1 z_2 z_3)^{-t-1} \left( \frac{ z_1^{t+2} }{ (z_1 - z_2) (z_1 - z_3) } + \frac{ z_2^{t+2} }{ (z_2 - z_1) (z_2 - z_3) } + \frac{ z_3^{t+2} }{ (z_3 - z_1) (z_3 - z_2) } \right)^3 d \z . \] 
We have to order the radii of the integration paths for each variable; we choose 
$ 0 < \epsilon_3 < \epsilon_2 < \epsilon_1 < 1 $. We heavily use this fact after 
multiplying out the cubic: integrating, for example, the term 
  \[ \frac{ z_1^{-t-1} z_2^{-t-1} z_3^{2t+5} }{ (z_3 - z_2)^3 (z_3 - z_1)^3 }  \] 
with respect to $z_3$ gives 0, as this function is analytic at the $z_3$-origin and $\left|z_1\right|, \left|z_2\right| > \epsilon_3$. 
After exploiting this observation for all the terms stemming from the cubic, 
the only integrals surviving are 
  \[ \frac{ 1 }{ (2 \pi i)^3 } \int \frac{ z_1^{2t+5} z_2^{-t-1} z_3^{-t-1} }{ (z_1 - z_2)^3 (z_1 - z_3)^3 } \ d \z \] 
and 
  \[ - \frac{ 3 }{ (2 \pi i)^3 } \int \frac{ z_1^{t+3} z_2 z_3^{-t-1} }{ (z_1 - z_2)^3 (z_1 - z_3)^2 (z_2 - z_3) } \ d \z \ . \] 
The first integral factors and yields, again by residue calculus, 
\begin{align*}
&\frac{ 1 }{ (2 \pi i)^3 } \int \frac{ z_1^{2t+5} z_2^{-t-1} z_3^{-t-1} }{ (z_1 - z_2)^3 (z_1 - z_3)^3 } \ d \z = \frac{ 1 }{ (2 \pi i)^3 } \int z_1^{2t+5} \left( \int \frac{ z^{-t-1} }{ (z_1 - z)^3 } \ dz \right)^2 dz_1 \\ 
&\qquad = \frac{ 1 }{ 2 \pi i } \int z_1^{2t+5} \left( - \frac 1 2 (-t-1)(-t-2) \, z_1^{-t-3} \right)^2 dz_1 = { \binom{ t+2 }{ 2 } }^2 \ . 
\end{align*} 
For the second integral, it is most efficient to integrate with respect to $z_2$ first: 
\begin{align*}
&- \frac{ 3 }{ (2 \pi i)^3 } \int \frac{ z_1^{t+3} z_2 z_3^{-t-1} }{ (z_1 - z_2)^3 (z_1 - z_3)^2 (z_2 - z_3) } \ d \z = - \frac{ 3 }{ (2 \pi i)^2 } \int \frac{ z_1^{t+3} z_3^{-t} }{ (z_1 - z_3)^5 } \ dz_3 \ dz_1 \\ 
&\qquad = - \frac{ 3 }{ 2 \pi i } \int z_1^{t+3} \, \frac 1 {4!} (-t) (-t-1) (-t-2) (-t-3) \, z_1^{-t-4} \ dz_1 = - 3 \binom{ t+3 }{ 4 } \ . 
\end{align*} 
Adding up the last two lines gives finally 
  \[ H_3(t) = { \binom{ t+2 }{ 2 } }^2 - 3 \binom{ t+3 }{ 4 } = \frac 1 8 t^4 + \frac 3 4 t^3 + \frac {15} 8 t^2 + \frac 9 4 t + 1 \ , \] 
which is equal to (\ref{h3}). 
To obtain the volume of $\B_3$, the leading term of $H_3$ has to be multiplied by the relative volume 
of the fundamental domain of the sublattice of $\Z^9$ in the affine space spanned by $\B_3$. This volume is 9; hence 
  \[ \vol \B_3 = \frac 9 8 \ . \] 
In general, it is not hard to prove (see, for example, the appendix of \cite{chanrobbinsbirkhoff}) 
that the relative volume of the fundamental domain of the sublattice of $\Z^{n^2}$ in the affine space spanned by $\B_n$ is $n^{n-1}$. 

The number of integrals we have to evaluate to compute $H_4$ is only slightly higher. By 
Theorem \ref{Hthm}, 
  \[ H_4 (t) = \frac 1 { ( 2 \pi i )^4 } \int_{ \left| z_1 \right| = \epsilon_1 } \int_{ \left| z_2 \right| = \epsilon_2 } \int_{ \left| z_3 \right| = \epsilon_3 } \int_{ \left| z_4 \right| = \epsilon_4 } ( z_1 z_2 z_3 z_4 )^{ -t-1 } \left( \sum_{k=1}^4 \frac{ z_k^{ t+3 } }{ \prod_{j \not= k} (z_k - z_j) } \right)^4 d \z \ . \] 
Again we have a choice of ordering the radii; we use $ 0 < \epsilon_4 <
\epsilon_3 < \epsilon_2 < \epsilon_1 < 1 $.  
After multiplying out the quartic, we have to calculate five integrals;
their evaluation---again  
straightforward by means of the residue theorem---is as follows. 
As before, we can `save' computation effort by choosing a particular order with
which we integrate and by factoring an integral if possible. 
\begin{align*}
& \frac{ 1 }{ (2 \pi i)^4 } \int \frac{ z_1^{3t+11} z_2^{-t-1} z_3^{-t-1}
z_4^{-t-1} }{ (z_1 - z_2)^4 (z_1 - z_3)^4 (z_1 - z_4)^4 } \ d \z = \frac{ 1
}{ (2 \pi i)^4 } \int z_1^{3t+11} \left( \int \frac{ z^{-t-1} }{ (z_1 - z)^4
} \ d z \right)^3 d z_1 \\
& \qquad = { \binom{ t+3 }{ 3 } }^3 ,
\displaybreak[0] \\
& - \frac{ 4 }{ (2 \pi i)^4 } \int \frac{ z_1^{2t+8} z_2^2 z_3^{-t-1}
z_4^{-t-1} }{ (z_1 - z_2)^4 (z_1 - z_3)^3 (z_1 - z_4)^3 (z_2 - z_3) (z_2 -
z_4)  } \ d \z \\ 
& \qquad = - \frac{ 4 }{ (2 \pi i)^4 } \int \frac{ z_1^{2t+8} z_2^2 }{ (z_1 -
z_2)^4 } \left( \int \frac{ z^{-t-1} }{ (z_1 - z)^3 (z_2 - z) } \ dz
\right)^2 dz_1 \ dz_2 \\  
& \qquad = \frac{ 4 }{ (2 \pi i)^2 } \int \frac{ z_1^{2t+8} z_2^{ -t+1 } }{
(z_1 - z_2)^7 } \left( 2 \binom{t+2}{2} \frac{ z_1^{ -t-3 } }{ z_1 - z_2 } +
2 (t+1) \frac{ z_1^{ -t-2 } }{ (z_1 - z_2)^2 } + 2 \frac{ z_1^{ -t-1 } }{
(z_1 - z_2)^3 } - \frac{ z_2^{ -t-1 } }{ (z_1 - z_2)^3 } \right) d \z \\  
& \qquad = 8 \binom{ t+2 }{ 2 } \binom{ t+5 }{ 7 } + 8 (t+1) \binom{ t+6 }{ 8
} + 8 \binom{ t+7 }{ 9 } - 4 \binom{ 2t+8 }{ 9 } \ ,\
\displaybreak[0] \\ 
& \frac{ 4 }{ (2 \pi i)^4 } \int \frac{ z_1^{2t+8} z_2^{-t-1} z_3^2
z_4^{-t-1} }{ (z_1 - z_2)^3 (z_1 - z_3)^4 (z_1 - z_4)^3 (z_2 - z_3) (z_3 -
z_4)  } \ d \z  \\  
& \qquad = \frac{ 4 }{ (2 \pi i)^3 } \int \frac{ z_1^{2t+8} z_2^{-t+1}
z_4^{-t-1} }{ (z_1 - z_2)^3 (z_1 - z_4)^7 (z_2 - z_4) } \ d \z \\  
& \qquad = \frac{ 4 }{ (2 \pi i)^2 } \int \frac{ z_1^{2t+8} z_4^{-t-1} }{
(z_1 - z_4)^7 } \left( \binom{t+2}{2} \frac{ z_1^{ -t-3 } }{ z_1 - z_4 } +
(t+1) \frac{ z_1^{ -t-2 } }{ (z_1 - z_4)^2 } + \frac{ z_1^{ -t-1 } }{ (z_1 -
z_4)^3 } \right) d \z \\  
& \qquad = 4 \left( \binom{ t+2 }{ 2 } \binom{ t+5 }{ 7 } + (t+1) \binom{ t+6
}{ 8 } + \binom{ t+7 }{ 9 } \right) \ , 
\displaybreak[0] \\ 
& \frac{ 6 }{ (2 \pi i)^4 } \int \frac{ z_1^{t+5} z_2^{t+5} z_3^{-t-1}
z_4^{-t-1} }{ (z_1 - z_2)^4 (z_1 - z_3)^2 (z_1 - z_4)^2 (z_2 - z_3)^2 (z_2 -
z_4)^2 } \ d \z  \\  
& \qquad = \frac{ 6 }{ (2 \pi i)^4 } \int \frac{ z_1^{t+5} z_2^{t+5} }{ (z_1
- z_2)^4 } \left( \int \frac{ z^{-t-1} }{ (z_1 - z)^2 (z_2 - z)^2 } \ dz
\right)^2 d z_1 \ d z_2 \\  
& \qquad = \frac{ 6 }{ (2 \pi i)^2 } \int \frac{ z_1^{t+5} z_2^{t+5} }{ (z_1
- z_2)^4 } \left( (t+1)^2 \frac{ z_2^{ -2t-4 } }{ (z_1 - z_2)^4 } - 4 (t+1)
\frac{ z_2^{ -2t-3 } }{ (z_1 - z_2)^5 } + 4 \frac{ z_2^{ -2t-2 } }{ (z_1 -
z_2)^6 } \right) d \z \\  
& \qquad = 6 (t+1)^2 \binom{ t+5 }{ 7 } - 24 (t+1) \binom{ t+5 }{ 8 } + 24
\binom{ t+5 }{ 9 } \ , 
\displaybreak[0] \\ 
& - \frac{ 12 }{ (2 \pi i)^4 } \int \frac{ z_1^{t+5} z_2^2 z_3^2 z_4^{-t-1}}
{(z_1 - z_2)^3 (z_1 - z_3)^3 (z_1 - z_4)^2 (z_2 - z_3)^2 (z_2 - z_4) (z_3 -
z_4) } \ d \z  \\ 
& \qquad = - \frac{ 12 }{ (2 \pi i)^3 } \int \frac{ z_1^{t+5} z_2^2
  z_4^{-t+1} }{ (z_1 - z_2)^3 (z_1 - z_4)^5 (z_2 - z_4)^3 } \ d \z \\  
& \qquad = - \frac{ 12 }{ (2 \pi i)^2 } \int \frac{ z_1^{t+5} z_4^{-t+1} }{
  (z_1 - z_4)^5 } \left( \frac 1 { (z_1 - z_4)^3 } + 6 \frac{ z_4 }{ (z_1 -
  z_4)^4 } + 6 \frac{ z_4^2 }{ (z_1 - z_4)^5 } \right) d \z \\
& \qquad = - 12 \binom{ t+5 }{7} - 72 \binom{ t+5 }{8} - 72 \binom{ t+5}{9} \ .
\end{align*} 
Adding them all up gives 
\[
H_4 (t) 
=  \frac {11}{11340} \,t^{9} + 
  \frac {11}{630} \,t^{8} + \frac {19}{135} \,t^{7} + \frac {2}{3} \,t^{6} +
  \frac {1109}{540} \,t^{5} + \frac {43}{10} \,t^{4} +  
  \frac {35117}{5670} \,t^{3} + \frac {379}{63} \,t^{2} + \frac {65}{18} \,t +
  1\\
\]
and hence 
  \[ \vol \B_4 = 4^3\cdot\frac {11}{11340}= \frac{176}{2835} \ . \] 


\section{Larger $n$ do} 

As we have seen in the examples, after multiplying out the integrand of
Theorem \ref{Hthm}, many of the terms do not contribute to the integral. The
following proposition will provide us with a general statement to that
effect.

For a rational function $f$ in $n$ variables $z_j$ we use the notation
$d_r(f)$ for the degree of $f$ in the variables $z_1,\dots,z_r$.

\begin{proposition}\label{int} 
Suppose $p_1,\dots,p_n$ are integers, $q_{jk}$ are nonnegative integers ($1
\leq j < k \leq n$), $1 > \epsilon_1 > \dots > \epsilon_n > 0$, and 
\[
f(z_1,\dots,z_n) = \frac{\prod_{1\leq j\leq n} z_j^{p_j}}{\prod_ {1 \leq j <
k \leq n} (z_j - z_k)^{q_{jk}}}\ .
\]
If $1\leq r\leq n$ and $d_r(f)<-r$ then
\[
\int_{ \left| z_1 \right| = \epsilon_1 } \cdots \int_{ \left| z_n \right|
= \epsilon_n } f(\z)\;d\z = 0\ .
\]
\end{proposition}

\emph{Proof}.  We need only show that
\[
\int_{ \left| z_1 \right| = \epsilon_1 } \cdots \int_{ \left| z_r \right|
= \epsilon_r } f(\z)\;d\z = 0\ ,\tag{$*$}
\]
for then we can integrate over all $n$ variables by first integrating over
$z_1\dots z_r$.

If $r=1$ then $f$, considered as a
function of $z_1$, is analytic outside the circle $\left| z_1 \right| =
\epsilon_1$ and has zero residue at infinity since its degree is less than
$-1$, so $(*)$ follows.

We continue by induction, so suppose $r>1$ and $(*)$ is true for
smaller $r$.  Suppose that $d_r(f)<-r$; we may assume $d_{r-1}(f)\ge -(r-1)$.
Note that $d_r(f)\ge d_{r-1}(f)+d$ where $d$ is the degree of $f$ in the
single variable $z_r$ (the discrepancy is the sum of the exponents $q_{jr}$
for $j<r$).  Hence $d\le d_r(f)-d_{r-1}(f)<-1$.  We consider $f$ as a
function of $z_r$ and apply the residue theorem to the region outside the
circle $\left| z_r \right| = \epsilon_r$.  As above $f$ has zero residue at
infinity, so we only need to consider the residues at the poles $z_j$ for
$j<r$.  Evaluating these residues converts the integral of $f$ into a
(possibly huge) linear combination of integrals of functions of the same
form as $f$, but in the $n-1$ variables
$z_1,\dots,z_{r-1},z_{r+1},\dots,z_n$.  If $g$ is any one of these functions
then we easily calculate $d_{r-1}(g)=d_r(f)+1 < -(r-1)$, and, by the
induction hypothesis, the integral of $g$ is zero.  \hfill {} $\Box$

From this proposition we obtain the starting point for our `algorithm.' 

\begin{corollary}\label{comp}
For $ 1 > \epsilon_1 > \dots > \epsilon_n > 0 $ and $t\ge0$, $H_n(t) =$
\[
\frac 1 { ( 2 \pi i )^{n} } \int_{ \left| z_1 \right| = \epsilon_1
} \cdots \int_{ \left| z_n \right| = \epsilon_n } ( z_1 \cdots z_{n} )^{
-t-1 } \sum_{ m_1 + \dots + m_n = n } \hspace{-.31in} \mbox{}^*
\hspace{.31in} \binom{ n }{ m_1, \dots , m_n } \prod_{ k=1 }^n \left( \frac{
z_k^{ t+n-1 } }{ \prod_{j \not= k} (z_k - z_j) } \right)^{ m_k } d \z \ ,
\]
where $\sum^*$ denotes that we only sum over those $n$-tuples of
non-negative integers satisfying $m_1 + \dots + m_n = n$ and $m_1 + \dots +
m_r > r$ if $1\le r<n$.
\end{corollary} 

\emph{Remark}.  The condition on $m_1,\dots,m_n$ can be visualized through
lattice paths from $(0,0)$ to $(n,n)$ using the steps $(1,m_1), (1,m_2),
\dots, (1,m_n)$.  The condition means that these paths stay strictly above the
diagonal (except at the start and end).\footnote{The authors thank Lou Billera for pointing out this
lattice-path interpretation.}

\emph{Proof}. By Theorem \ref{Hthm}, 
\begin{align*}
H_n (t) &= \frac 1 { ( 2 \pi i )^{n} } \int \cdots \int ( z_1 \cdots
z_{n} )^{ -t-1 } \left( \sum_{k=1}^n \frac{ z_k^{ t+n-1 } }{ \prod_{j \not=
k} (z_k - z_j) } \right)^n d \z \\
&= \frac 1 { ( 2 \pi i )^{n} } \int \cdots \int ( z_1 \cdots z_{n} )^{ -t-1 } \sum_{ m_1 + \dots + m_n = n }
\binom{ n }{ m_1, \dots , m_n } \prod_{ k=1 }^n \left( \frac{ z_k^{ t+n-1 }
}{ \prod_{j \not= k} (z_k - z_j) } \right)^{ m_k } d \z \ . 
\end{align*}
We select a partition $m_1+\dots+m_n=n$ and rewrite the corresponding
integrand in the language of Proposition \ref{int}:
\[
p_j = (n-1)m_j + t(m_j-1) - 1\ ,\qquad q_{jk} = m_j+m_k\ .
\]
Now suppose $1\le r<n$. The degree of the denominator in $z_1,\dots,z_r$ is
\begin{align*} 
&\sum_{j=1}^r\sum_{k=j+1}^n (m_j+m_k) =
(n-1)\sum_{j=1}^r m_j + r\sum_{j=r+1}^nm_j = (n-1)\sum_{j=1}^r m_j +nr -
r\sum_{j=1}^r m_j \\ 
&\qquad = (n-1)\sum_{j=1}^r m_j +nr -r\sum_{j=1}^r(m_j-1) -r^2 \ . 
\end{align*}
We subtract this from $\sum_{j=1}^r p_j$ to get
\[
d_r(f) = (t+r)\sum_{j=1}^r(m_j-1) - r(n-r) - r\ .
\]
Here $t$ is non-negative and $r(n-r)$ is positive, so Proposition \ref{int}
implies that the integral is zero unless $\sum_{j=1}^r(m_j-1)>0$.\hfill {}
$\Box$

Theoretically, Corollary \ref{comp} tells us what we have to do to compute
$H_n$.  Thus the integrals we computed in our calculations of $H_3$ and
$H_4$ are exactly the non-zero integrals according to this result.
For practical purposes, however, the statement is almost worthless for larger values 
of $n$.  The first problem is that the number of terms in 
the sum equals the $(n-1)^{ \text{th} }$ Catalan number 
  \[ \frac 1 n \binom{ 2(n-1) }{ n-1 } = \frac{ (2n-2)! }{ n! (n-1)! } \ , \] 
which grows exponentially with $n$. 
Another slippery point is the evaluation of each integral. As we have seen in the 
examples, and as can be easily seen for the general case, we can compute each integral 
step by step one variable at a time. However, this means at each step we convert a 
rational function into a sum of rational functions (of one variable less) by means of 
the residue theorem. Again, this means that the number of (single-variable) integrals we have to compute 
grows immensely as $n$ increases. 
In fact, if we just `feed' the statement of Corollary \ref{comp} into a computer and 
tell it to integrate each summand, say, starting with $z_1$, then $z_2$, and so on, 
the computation time explodes once one tries $n=7$ or $8$. 
We feel that computationally this is as involved as calculating a sufficient number of 
values of $H_n$ and then interpolating this polynomial. 
However, complex analysis allows us some shortcuts which turn out to speed up the computation 
by a huge factor and which make Corollary \ref{comp} valuable, even from a computational 
perspective. These `tricks' all showed up already in the examples and include 
\vspace{-4mm} 
  \begin{enumerate} 
    \item realizing when a function is analytic at the $z_k$-origin, 
    \item trying to choose the most efficient order of variables to
    integrate (based on estimating how many terms will be generated
    by the residue calculations, for each available variable), and
    \item\label{last} factoring the integral if some of the variables appear in a symmetric fashion. 
  \end{enumerate} 
\vspace{-4mm} 
Finally, if we are only interested in the volume of $\B_n$, we may also be 
\vspace{-4mm} 
  \begin{enumerate} 
    \item[4.] suppressing a particular integral if it does not contribute to the leading term of $H_n$. 
  \end{enumerate} 
It is worth noting that each of these computational `speed-ups' decreases the total computation time substantially. 
By applying them to Corollary \ref{comp}, we implemented a 
C{\raisebox{2pt}{\scriptsize$++$} program 
for the specific functions we have to integrate to compute $H_n$.  
We were able to verify all previously known polynomials ($n \leq 8$) and to compute $\vol \B_9$, $\vol \B_{10}$, 
and $H_9$. 
The results of our calculations,
including the polynomials $H_n$ for $n\le9$, and the source code for our
program are available at {\tt www.math.binghamton.edu/dennis/Birkhoff}.
The following table gives the volumes together with the respective computing time (on a 1GHz PC running under Linux). 
Note that computing the full polynomial $H_n$ takes longer, as we cannot make use of shortcut \#4. 
In fact, the computation of $H_9$ took about 325 days of computer
time, although the 
elapsed time was only about two weeks since we distributed the calculations
among a number of machines (8 -- 40, depending on availability).

\begin{center}
\begin{tabular}{c|c|c} $n$ & $\vol \B_n$ & time \\ \hline 
1 & 1 &$<.01$ sec \\ \hline 
2 & 2 &$<.01$ sec\\ \hline 
3 & \Large $ \rule[-10pt]{0pt}{25pt}\frac 9 8$ &$<.01$ sec\\ \hline 
4 & \Large $ \rule[-10pt]{0pt}{25pt}\frac{176}{2835} $ &$<.01$ sec\\ \hline 
5 & \Large $ \rule[-10pt]{0pt}{25pt}\frac{23590375}{167382319104} $ &$<.01$
sec\\ \hline  
6 & \Large $ \rule[-10pt]{0pt}{25pt}\frac{9700106723}{1319281996032 \cdot
10^6} $ &$.18$ sec\\ \hline 
7 & \Large $
\rule[-10pt]{0pt}{25pt}\frac{77436678274508929033}{137302963682235238399868928
\cdot 10^8} $ &$15$ sec\\ \hline 
8 & \Large $
\rule[-10pt]{0pt}{25pt}\frac{5562533838576105333259507434329}{12589036260095477950081480942693339803308928
\cdot 10^{10}} $ &$54$ min\\ \hline 
9 & \Large $
\rule[-10pt]{0pt}{25pt}\frac{559498129702796022246895686372766052475496691}{92692623409952636498965146712806984296051951329202419606108477153345536
\cdot 10^{14}} $ &$317$ hr
\end{tabular}
\end{center} 

After this article was submitted for publication, we used the idle time on our 
departmental machines to compute the volume of $\B_{10}$. After a computation 
time of 6160 days, or almost 17 years (again scaled to a 1GHz processor), we 
obtained 

$\vol \B_{10} = $ 
\scriptsize 
\[ \frac{727291284016786420977508457990121862548823260052557333386607889} 
{8281608601067668551256763187968727293446224635330894226779807213
88055739956270293750883504892820848640000000} \ . 
\] 
\normalsize 

We hope that the reader will take this immense computing time as a
challenge to improve our algorithms, and work towards $n=11$.


\section{An outlook towards transportation polytopes} 

The Birkhoff polytopes are special cases of transportation polytopes, which are defined below. 
The study of this class of polytopes, which are naturally at least as fascinating as the Birkhoff polytopes, 
was motivated by problems in linear programming; for combinatorial properties see, for example, \cite{diaconisgangoli,kleewitzgall}. 
The goal of this section is to show how our methods can be applied in this more general setting. 

Fix positive real numbers $a_1,\dots,a_m,b_1,\dots,b_n$ such that $a_1+\dots+a_m = b_1+\dots+b_n$. 
Let $\a=(a_1,\dots,a_m), \b=(b_1,\dots,b_n)$. 
  \[ \T_{\a,\b} = \left\{ \left( \begin{array}{ccc} x_{11} & \cdots & x_{1n} \\ \vdots & & \vdots \\ x_{m1} & \dots & x_{mn} \end{array} \right) \in \R^{mn} : \ x_{jk} \geq 0 , \begin{array}{l} \sum_k x_{jk} = a_j \text{ for all } 1 \leq j \leq m , \\ \sum_j x_{jk} = b_k \text{ for all } 1 \leq k \leq n \end{array} \right\} \] 
is the set of \emph{solutions to the transportation problem with parameters} $\a,\b$. 
It is a convex polytope of dimension $(m-1)(n-1)$ in $\R^{mn}$---hence we
refer to $\T_{\a,\b}$ as a \emph{transportation polytopes}. 
Note that $\B_n = \T_{\1,\1}$ where $\1 = (1,\dots,1) \in \Z^n$. 
We will be exclusively interested in the case when $a_1,\dots,a_m,b_1,\dots,b_n$ are integers. 
One reason for this is that the vertices of $\T_{\a,\b}$ are then integral. 
Let 
  \[ T(\a,\b) = T(a_1,\dots,a_m,b_1,\dots,b_n) = \# \left( \T_{\a,\b} \cap \Z^{mn} \right) \] 
denote the number of integer points in the transportation polytope with integral parameters $\a,\b$. 
Geometrically, each of these parameters determines the position of a hyperplane bounding the polytope $ \T_{\a,\b} $. 
It is well known that $ T(\a,\b) $ is a piecewise-defined polynomial in $a_1,\dots,a_m,b_1,\dots,b_n$. 
The regions in which $ T(\a,\b) $ is a polynomial depends on the normal cones of the polytopes 
involved; for details, we refer to \cite{sturmfelsvectorpartition}. 

Moreover, one can derive results for the counting function $ T(\a,\b) $ which are 
`higher-dimensional' extensions of {\rm (\ref{reczeros})}, by application of the 
main theorem in \cite{evencloser}, which in turn is a generalization of the 
Ehrhart-Macdonald reciprocity theorem (\ref{rec}). To state this theorem, we need to 
define the following integer-point counting function. Suppose $\P$ is a rational polytope 
given in the form (\ref{polytope}); denote by 
$\P(\t) = \left\{ \x \in \R^d : \ \A \x = \t \right\} $ 
a polytope which we obtain from $\P = \P(\b)$ by translating (some of) its bounding hyperplanes. 
(Classical Ehrhart dilation is the special case $\t = t \b$.) 
Let 
  \[ \L_\P (\t) = \# \left( \P(\t) \cap \Z^d \right) \ . \] 
\cite[Theorem 4]{evencloser} states that $\L_\P$ and $\L_{\P^\circ}$ are piecewise-defined 
(multivariable) quasipolynomials satisfying 
  \[ \L_\P(-\t) = (-1)^{ \dim \P } \L_{\P^\circ}(\t) \ . \] 
As easily as (\ref{reczeros}) follows from (\ref{rec}), this reciprocity theorem 
yields a symmetry result for the transportation counting function. Let
$\1_d$ denote the $d$-dimensional 
vector all of whose entries are one. Then  
  \[ T(-\a-n\1_m,-\b-m\1_n) = (-1)^{ mn - 1 } \ T(\a,\b) \ . \] 

We finally turn to the problem of writing $T(\a,\b)$ in form of an integral. 
As with $\B_n$, we can view $\T_{\a,\b}$ as given in the form of (\ref{polytope}) and apply the philosophy of 
Theorem \ref{intthm} to obtain 
  \[ T(\a,\b) = \frac 1 { ( 2 \pi i )^{m+n} } \int\cdots\int \frac{ z_1^{ -a_1 - 1 } \cdots z_m^{ -a_m - 1 } w_1^{ -b_1 - 1 } \cdots w_n^{ -b_n - 1 } }{ \prod_{ { 1 \leq j \leq m } \atop { 1 \leq k \leq n } } \left( 1 - z_j w_k \right) } \ d \z \ d \w \ . \] 
Again the integral with respect to each (one-dimensional) variable is over a circle with radius $<1$ centered at 0, 
and all appearing radii are different. As in the Birkhoff case, we can separate, say, the $\w$-variables and obtain 
  \[ T(\a,\b) = \frac 1 { ( 2 \pi i )^{m+n} } \int\cdots\int z_1^{ -a_1 - 1 } \cdots z_m^{ -a_m - 1 } \prod_{ k=1 }^n \int \frac{ w_k^{ -b_k - 1 } }{ \prod_{ j=1 }^m \left( 1 - z_j w_k \right) } \ d w_k \ d \z \ . \] 
And as with the Birkhoff polytope, the innermost integrals are easy to compute by means of the residue theorem. 
This yields the following statement which generalizes Theorem \ref{Hthm}. 
\begin{theorem}\label{Tthm} For any distinct $ 0 < \epsilon_1, \dots, \epsilon_m < 1 $, 
  \[ T(\a,\b) = \frac 1 { ( 2 \pi i )^{m} } \int_{ \left| z_1 \right| = \epsilon_1 } \cdots \int_{ \left| z_m \right| = \epsilon_m } z_1^{ -a_1 - 1 } \cdots z_m^{ -a_m - 1 } \prod_{ k=1 }^n \sum_{j=1}^m \frac{ z_j^{ b_k+m-1 } }{ \prod_{l \not= j} (z_j - z_l) } \ d \z \ . \] 
\end{theorem} 

\emph{Remark}. As with Theorem \ref{Hthm}, it can be proved from the form of the integrand 
that $ T(\a,\b) $ is indeed a piecewise-defined polynomial in $a_1,\dots,a_m,b_1,\dots,b_n$. 

Theoretically we could now use this theorem to produce formulas for $T(\a,\b)$ just as we did for $H_n(t)$. 
There is one major difference: $T(\a,\b)$ is only a \emph{piecewise-defined} polynomial. In fact, we 
can see this from the form of the integral in Theorem \ref{Tthm}: whether we will get a nonzero contribution 
at a certain step in the computation depends heavily on the relationship between $a_1,\dots,a_m,b_1,\dots,b_n$. 

The fact that the counting function $T(\a,\b)$ is of a somewhat more delicate nature 
naturally has computational consequences. We feel that providing any general results on this function 
would go beyond the scope of this article and will hopefully be the subject of a future project. 
On the other hand, we adjusted our algorithm to compute values of $T(\a,\b)$ for three (fixed) pairs $(\a,\b)$, 
which have been previously computed by Mount \cite{mount} and DeLoera and Sturmfels \cite{deloerasturmfels}. 
The first example is 
\begin{align*} 
&T((3046,5173,6116,10928), (182,778,3635,9558,11110)) = \\ 
&\qquad 23196436596128897574829611531938753 
\end{align*} 
The authors reported their computation took 20 minutes (Mount) or 10
minutes (DeLoera--Sturmfels) \cite{deloerasturmfels,mount}. 
We computed this number in 0.2 seconds, based on Theorem \ref{Tthm}.
A similar phenomenon happens with 
\begin{align*} 
&T((338106,574203,678876,1213008), (20202,142746,410755,1007773,1222717)) = \\ 
&\qquad 316052820930116909459822049052149787748004963058022997262397 
\end{align*} 
and 
\begin{align*} 
&T((30201,59791,70017,41731,58270), (81016,68993,47000,43001,20000)) = \\ 
&\qquad 24640538268151981086397018033422264050757251133401758112509495633028 
\end{align*} 
which reportedly took 35 minutes/10 days with the DeLoera--Sturmfels
algorithm \cite{deloerasturmfels}, 0.3/2.9 seconds with ours.  These timing
comparisons ignore differences in machine speeds and implementation of
the algorithms, but suggest that our methods are considerably more efficient.

\vspace{1cm}
{\bf Acknowledgements}. We would like to thank the referees for carefully reading through our paper and software, 
and for their many helpful comments. 


\newpage 
\bibliographystyle{plain}

\providecommand{\bysame}{\leavevmode\hbox to3em{\hrulefill}\thinspace}
\providecommand{\MR}{\relax\ifhmode\unskip\space\fi MR }
\providecommand{\MRhref}[2]{%
  \href{http://www.ams.org/mathscinet-getitem?mr=#1}{#2}
}
\providecommand{\href}[2]{#2}

\sc Department of Mathematical Sciences\\
    State University of New York\\
    Binghamton, NY 13902-6000\\
\tt matthias@math.binghamton.edu \\ 
    dennis@math.binghamton.edu \\ 

\end{document}